\documentclass[12pt]{article}

\title{Some applications of duality\\
for L\'evy processes in a half-line} 
\author{Jean Bertoin\thanks{Laboratoire de Probabilit\'es, 
UPMC, 175 rue du Chevaleret, 75013 Paris; and DMA, ENS Paris, France. Email:
jean.bertoin@upmc.fr} \and 
Mladen Savov
\thanks{ Department of Statistics,
1 South Parks Road,
Oxford, OX1 3TG
United Kingdom.
Email: savov@stats.ox.ac.uk} 
 }

\date{}

\usepackage{amsfonts,amsmath,amssymb,bbm}
\usepackage{setspace} 
\def\proof{\noindent{\bf Proof:}\hskip10pt}        
\def\QED{\hfill $\Box$}

\font\tenmath=msbm10 scaled 1200
\font\sevenmath=msbm7 scaled 1200
\font\Phiivemath=msbm5 scaled 1200
\newfam\mathfam \textfont\mathfam=\tenmath
\scriptfont\mathfam=\sevenmath \scriptscriptfont\mathfam=\Phiivemath

\def \\ { \cr }
\def\R{\mathbb{R}}
\def \1{1 \mkern -6mu 1}

\def\P{\mathbb{P}}

\def\D{\mathbb{D}}

\def \up{^{\uparrow}}
\def \e{{\rm e}}
\def \d{{\rm d}}

\newtheorem{theorem}{Theorem}

\newtheorem{lemma}{Lemma}
\newtheorem{corollary}{Corollary}
\begin{document}

\maketitle

\begin{abstract}
The central result of this paper is an analytic duality relation for real-valued  L\'evy processes killed upon exiting a half-line. By Nagasawa's theorem, this yields a remarkable time-reversal identity
involving the L\'evy process conditioned to stay positive. As examples of applications, we construct a version of the L\'evy process indexed by the entire real line and started from $-\infty$
which enjoys a natural spatial-stationarity property, and point out that the latter leads to a natural Lamperti-type representation for self-similar Markov processes in $(0,\infty)$ started from the entrance point $0+$.
  \end{abstract}

\begin{section}{Introduction}
A celebrated result due to David Williams (cf. Theorem 3.4 in \cite{DW}) can be stated as follows. Consider a real Brownian motion $(B^x_t)_{t\geq 0}$ started from some level $x>0$ and $T=\inf\{t\geq 0: B^x_t\leq 0\}$ its first exit-time from $(0,\infty)$. Then the process
$(B^x_{T-t})_{0\leq t < T}$ obtained by time-reversing the Brownian path at time $T$, has the same distribution as a three-dimensional Bessel process started from $0$ and killed at the time of its last-passage at level $x$. This relation should be viewed as the probabilistic counterpart of an analytic duality between the transition probabilities $q_t(x,\d y)$ of the Brownian motion killed upon exiting $(0,\infty)$ and $p_t\up(x,\d y)$ of the three-dimensional Bessel process. Specifically, there is the identity
\begin{equation}\label{eq1}
 q_t (x,\d y) x\d x= p\up_t(y,\d x) y \d y \,,\qquad x,y\in(0,\infty).
 \end{equation}
One further observes that the duality measure $x\d x$ on $(0,\infty)$ coincides with 
the potential measure $U\up(0,\d x)=\int_0^{\infty}\d t p\up_t(0,\d x)$ of the three-dimensional Bessel process started from $0$, and the time-reversal identity of Williams then follows from a general result of Nagasawa \cite{Na} for Markov processes in duality. We also refer to 
Az\'ema \cite{Az} and Chung and Walsh \cite{CW} for further seminal contributions in this area. 

It has been observed in Theorem VII.18 in \cite{Be} that similar arguments
can be applied  to L\'evy processes with no positive jumps, and yield an extension of
 Williams' time-reversal identity in that setting. More precisely, consider a L\'evy process $\xi$ with no positive jumps; the role of the Brownian motion is now played by $\hat \xi=-\xi$, and that of the three-dimensional Bessel process by $\xi\up$, which should be thought of as $\xi$ conditioned to stay positive (in general such a conditioning is singular and has to be understood in terms of Doob's $h$-transform). Of course, the absence of positive jumps of $\xi$ is crucial as it ensures that the downwards passages for $\hat\xi$  occur continuously. The central result of the present work is that a duality identity extending \eqref{eq1} holds for general L\'evy processes (possibly with positive jumps), and as a consequence so does the remarkable time-reversal identity for L\'evy processes which do not tend to $-\infty$. A fundamental feature of this extension is the possibility of downwards crossings by a jump for $\hat \xi$, so in general
 the L\'evy process $\hat\xi$ and the version of $\xi$ conditioned to stay positive, $\xi\up$, have to start  from appropriate random locations in $[0,\infty)$.

This duality relation has a number of applications, some of which have already been observed in the literature. In particular, it enables us to construct a process $(\xi_t)_{t\in\R}$ 
indexed by the real line that fulfills a spatial invariance property and which may be thought of as a version of the L\'evy process $\xi$ started from $-\infty$. More precisely,  it  appears as the limit in distribution as $x\to -\infty$ of the L\'evy process started at time $0$ from $x$ and shifted in time at the instant of its first entrance in $(0,\infty)$. In this vein, we point at a remarkable representation of positive self-similar Markov processes $(X_t)_{t\geq 0}$ in $(0,\infty)$ started from the boundary point $0+$ as a time-change of $\exp (\xi)$, which extends the classical construction due to Lamperti \cite{La} when $X$ starts from a strictly positive position.

This paper is organized as follows. In the next section we first recall some useful  elements of fluctuation theory for L\'evy processes, and then present the key duality relation. After discussing the classical weak convergence of the under and over shoots in the framework of renewal theory applied to the ladder height process, we finally use Nagasawa's theorem to establish 
an identity involving time-reversed processes which provides the probabilistic counterpart of the duality relation. Section 3 is devoted to applications to limit theorems. We first observe that there is a natural version of the L\'evy process indexed by the entire real line which enjoys a remarkable spatial-stationarity property. Then we show that this process appears as the limit as $z\to -\infty$ of the genuine L\'evy process started from $z$ and shifted in time at its first entrance in $(0,\infty)$. This weak limit theorem
encompasses the classical convergence of the under and over shoots, and the combination with Lamperti's transformation points at
a simple approach for studying the entrance boundary of positive self-similar Markov processes.

 \end{section}

 \begin{section}{Duality and time-reversal in a half-line}
 \subsection{Some notation and preliminaries}
 We introduce some background on L\'evy processes and fluctuation theory that will be needed here, referring the reader to \cite{Be}, \cite{Do} or \cite{Ky} for a complete account. {\bf We implicitly exclude the compound Poisson processes} (merely to avoid discussing periodicity).
 
 We shall use the canonical notation :  the probability space is chosen to be   $\Omega=\D([0,\infty),\R)$, the space of c\`adl\`ag paths  endowed with the Borel sigma-field generated by Skorohod's topology, and $\xi=(\xi_t)_{t\geq 0}$
 is the coordinate process, i.e. $\xi_t(\omega)=\omega(t)$.
 Our building block is a probability measure $P$ on $\Omega$ for which $\xi$ is a L\'evy process, i.e. $\xi$ has independent and homogeneous increments and starts from $\xi_0=0$ a.s.  We  write $\Pi$ for the L\'evy measure, which specifies the intensity of the jumps.  We also denote  by $\hat P$ the image of $P$ by the map $\omega\to \hat \omega=-\omega$. In other words, $\hat P$ is the law of the dual L\'evy process $\hat \xi=-\xi$ under $P$.
 
 Killing the paths at their first-exit time from the upper half-line
  $$T=\inf\{t\geq 0: \xi_t\leq 0\}$$
yields two sub-Markovian transition probabilities on $(0,\infty)$
$$p_t(x,\d y)=P_x(\xi_t\in \d y, t<T)\ \hbox{and}\ \hat p_t(x,\d y)=\hat P_x(\xi_t\in \d y, t<T)\,,$$
where $P_x$ and $\hat P_x$ denote the law of $x+\xi$ under $P$ and under $\hat P$, respectively. 
We write
 $$U(x,\d y)=\int_0^{\infty}\d t p_t(x,\d y)$$
for the potential measure of the L\'evy process killed when exiting $(0,\infty)$.

 Recall that under $P$,  the reflected process 
 $(\sup_{0\leq s\leq t}\xi_s-\xi_t)_{t\geq 0}$ is a Feller process in $[0,\infty)$  which possesses a local time $(L_t)_{t\geq0}$ at level $0$.  The (ascending) ladder time is defined as  the right-continuous inverse of $L$, viz. $L^{-1}(t)=\inf\{s\geq 0: L_s>t\}$
 and the ladder height process 
 $H_+$ by 
 $$H_+(t)=\xi_{L^{-1}(t)}=\sup_{0\leq s\leq L^{-1}(t)}\xi_s\,, \qquad \hbox{ whenever } L^{-1}(t)<\infty\,.$$ 
Here, we use the convention $\inf\varnothing =\infty$ and $H_+(t)=\infty$ when $L_{\infty}\leq t$.
 It is well-known that $H_+$ is a subordinator (killed at time $L_{\infty}$ when the L\'evy process tends to $-\infty$).  We denote its drift coefficient by $a_+\geq 0$ and its L\'evy measure by $\mu_+$, so that for every $q,t\geq 0$, 
 $$E\left(\exp(-qH_+(t))\right) = 
 \exp\left(-t \left(a_+q + \int_{(0,\infty]}\mu_+(\d x) (1-\e^{-qx})\right)\right)\,.$$
 Here we agree that $\exp(-qH_+(t))=0$ when $L^{-1}(t)=\infty$, and $\mu_+(\{\infty\})$
 corresponds to the killing rate of $H_+$.
We write $U_+$ for its renewal function, viz.
 $$U_+(x)=\int_0^{\infty}P(H_+(t)\leq x, L^{-1}(t)<\infty)\d t\,,\qquad x\in [0,\infty)\,.$$
 
 We also consider the dual ladder $H_-$, that is the ladder height of the dual L\'evy process $\hat \xi$ and denote by  $U_-$ its renewal function. 
 According to Silverstein \cite{Si} (see also Theorem VI.20 in \cite{Be}), there is the 
 remarkable identity
\begin{equation} \label{eq2}
U(x,\d y)= \int_{(x-y)^+}^x U_-(\d z) U_+(\d y+z-x)\,.
\end{equation}
More precisely, Silverstein's identity often appears with an additional constant factor $c>0$ 
in the right-hand side of \eqref{eq2}, which depends on the choice of the normalization that has been used to define the local times at $0$ of the reflected L\'evy processes. 
We thus implicitly assume that the local times have been normalized so that \eqref{eq2} holds exactly.

Silverstein \cite{Si} also
observed that the renewal function $U_-$ is harmonic for the semigroup induced by $(p_t)_{t\geq 0}$, i.e.
$$U_-(x)=\int_0^{\infty} p_t(x,\d y) U_-(y)\,,\qquad x>0\,.$$
Following Doob, this enables us to construct  (conservative) Markovian transition functions
$$p\up_t(x,\d y) = \frac{U_-(y)}{U_-(x)}p_t(x,\d y)\,.$$
The distribution of the Markov process on $(0,\infty)$ started from $x>0$
and with  transition functions $(p\up_t)_{t\geq 0}$ 
will be denoted by $P\up_x$;  roughly speaking, $P\up_x$ should be viewed as the law of the L\'evy process started from $x$ and conditioned to stay positive. 
Indeed, if the L\'evy process tends to $\infty$, that is $P_x(T=\infty)>0$ for some (and then all)
$x>0$, then it is easily seen that there exists some constant $c'>0$ such that 
$U_-(x) = c' P_x(T=\infty)$ and hence $p\up_t(x,\d y) $ coincides with  the transition probability
of the L\'evy process conditioned to stay positive in the usual sense. 
Finally, we denote the potential measure of the L\'evy process started from $x$ and conditioned to stay positive by
\begin{equation}\label{eq3}
U\up(x,\d y) = \int_0^{\infty}\d t p\up_t(x,\d y) = \frac{U_-(y)}{U_-(x)} U(x,\d y)\,. 
\end{equation}
 
When
 $0$ is regular upwards, in the sense that 
$P(\sup_{0\leq s\leq \varepsilon}\xi_s>0)=1$ for any $\varepsilon>0$,
it is known  from a work of Chaumont and Doney \cite{CD}
that $P_x\up$ has a weak limit as $x\to 0+$, which we denote by $P_0\up$.
More precisely $P_0\up(\xi_t>0)=1$ for all $t>0$ and 
under $P_0\up$ the canonical process $\xi$ remains Markovian (as a matter of fact, Fellerian) with transition probabilities $(p\up_t)_{t\geq 0}$.  We shall need the following result when the ladder height subordinator $H_+$ has a strictly positive drift coefficient  (we refer to Vigon \cite{Vi} for an explicit  necessary and sufficient condition for this to happen).

\begin{lemma}\label{L1} Suppose that $a_+>0$. Then the following holds :

\noindent{\rm (i)} The renewal function $U_+$ has a continuous derivative $u_+$ which is strictly positive everywhere with $u_+(0)=1/a_+$. For every $x>0$, the probability that the dual L\'evy process started from $x$ exits $(0,\infty)$ for the first time continuously is
$$\hat P_x(\xi_T=\xi_{T-}, T<\infty)= a_+u_+(x)\,.$$
Further, this quantity converges to $a_+/E(H_+(1))$ when $x\to\infty$.

\noindent{\rm (ii)} There is the identity
$$U\up(0,\d y)=U_-(y)u_+(y)\d y\,.$$
\end{lemma}

\proof  The probability under $\hat P_x$ that  $\xi$ is continuous at its  first exit time from $(0,\infty)$  coincides with the probability that the ladder height subordinator $H_+$ crosses the level $x$ continuously.  The first assertion in (i) merely rephrases a result of Neveu which is stated as Theorem III.5 in \cite{Be}, while the last one is a consequence of the renewal theorem for subordinators (see, e.g. Proposition 3.3 in \cite{BeSF}). 

So we focus on (ii).  The existence of a regular density $u_+$ for the renewal function $U_+$ enables us to rewrite \eqref{eq2} in the form
 $$ U(x,\d y)/\d y= \int_{(x-y)^+}^x U_-(\d z) u_+(y+z-x)\,.$$
Combining with \eqref{eq3} readily yields the desired formula. \QED
 
 \subsection{A duality relation}
 In order to state the duality relation that lies at the heart of this work, we still need 
 one more notation. We introduce the measure
\begin{equation}\label{eq4}
m(\d x)=\overline\Pi(x)U_-(x) \d x + a_+\delta_0(\d x)\,,\qquad x\in[0,\infty)
\end{equation}
 where $\overline\Pi(x)=\Pi((x,\infty))$ is the (upper) tail distribution of the L\'evy measure $\Pi$
 and  $a_+$ the drift coefficient of the ladder height subordinator $H_+$. 
 
 \begin{theorem}\label{T1}
 {\rm (i)} There is the duality identity
 $$p\up_t(x,\d y) U_-(x)\d x = \hat p_t(y, \d x)U_-(y) \d y\,.$$
 
 \noindent {\rm (ii)}  The duality measure $U_-(x)\d x $ can be expressed as
 $$U_-(x)\d x = \int_{[0,\infty)}m(\d y)U\up(y,\d x)\,.$$
  \end{theorem}
 
 \proof
 (i) This follows immediately from Hunt's switching identity  
 $$p_t(x, \d y) \d x= \hat p_t(y, \d x)\d y$$
 (see Theorem II.5 in \cite{Be}) 
 and the definition of $p\up_t(x,\d y)$. Note that this has already been pointed out in the proof of 
 Theorem 4 of Chaumont \cite{Ch}.
 
 (ii) From  \eqref{eq3}  we get
 $$\int_{[0,\infty)}\d y \overline\Pi(y)U_-(y)U\up(y,\d x)= 
 U_-(x)\int_{[0,\infty)}\d y \overline\Pi(y)U(y,\d x).$$
 On the other hand, Hunt's switching identity gives
 $$\int_{[0,\infty)}\d y \overline\Pi(y)U(y,\d x) = \left( \int_{[0,\infty)}\hat U(x,\d y) \overline\Pi(y)
 \right) \d x$$
 where 
 $$\hat U(x,\d y)=\int_0^{\infty} \d t \, \hat p_t(x, \d y)$$
 is the potential measure of the dual L\'evy process $\hat \xi =-\xi$ killed upon exiting $(0,\infty)$. 
 
 A standard argument based on the Poissonian structure of the jumps of the dual  L\'evy process 
 and the compensation formula for Poisson point processes shows that 
 $$\int_{[0,\infty)}\hat U(x,\d y) \overline\Pi(y) = \hat P_x(\xi \hbox{ exits from $(0,\infty)$ by a jump})\,.$$ 
 See for instance the proof of Proposition III.2 in \cite{Be}. 
When the drift coefficient $a_+$ of the ladder height subordinator $H_+$ (under $P$) is zero, the probability above is one according to a celebrated result due to Kesten which is stated as Theorem II.4 in \cite{Be}, and the proof is complete.
When $a_+>0$, we deduce from Lemma \ref{L1} that
$$a_+U\up(0,\d x)/\d x= U_-(x) \hat P_x(\xi \hbox{ exits from $(0,\infty)$ continuously })\,,$$
which yields the conclusion. \QED
 
 \subsection{Weak convergence of the over and under shoots}
 The probabilistic interpretation of the duality identity (Theorem \ref{T1}) requires the measure $m$ defined by \eqref{eq4} to be finite. The following claims are due to Vigon  \cite{Vi} (see (5.3.4) in  \cite{Do}) and 
 Doney and Maller (see Theorem 8 in \cite{DM}), respectively. 
 
 \begin{lemma} \label{L2} The mass of the measure $m$ coincides with the mean ladder height, i.e. 
 $$m([0,\infty))= E(H_+(1))\,.$$
 This quantity is  finite if and only if $\xi_1\in L^1(P)$
 and either $E(\xi_1)>0$  or  $E(\xi_1)=0$ and 
$$\int_{[1,\infty)}\d x \frac{x \overline \Pi(x)}{\int_0^x\d y \int_y^{\infty}\d z \Pi((-\infty,-z))}<\infty\,.$$
 \end{lemma}
 
 {\bf We assume that $E(H_+(1))<\infty$ throughout the rest of this work, except at the beginning of Section 3.2}. We stress that this rules out the case when the L\'evy process tends to $-\infty$; in particular the ascending ladder processes are not defective. 
 
 Next, we introduce the probability measure $\rho$ on $[0,\infty)^2$ as
 \begin{equation}\label{eq5}
 \rho(\d x, \d y)= \frac{1}{E(H_+(1))} (U_-(x)\Pi(x+\d y) \d x+ a_+\delta_0(\d x)\delta_0(\d y))\,,
 \end{equation}
 and  write $\rho_1$ and $\rho_2$ for the marginal laws of $\rho$ :
\begin{eqnarray*}
\rho_1(\d x)&=& \frac{1}{E(H_+(1))} \left(U_-(x)\overline \Pi(x) \d x + a_+\delta_0(\d x)\right)\,, \\
 \rho_2(\d y)&=& \frac{1}{E(H_+(1))} \left( a_+\delta_0(\d y) +\int_0^{\infty}\d x U_-(x) \Pi(x+ \d y)\right) \,.
\end{eqnarray*}
Note that the first marginal $\rho_1$ coincides with the measure $m$ in \eqref{eq4} normalized to be a probability, and that an integration by parts gives
$$\int_0^{\infty}\d x U_-(x) \Pi(x+ \d y) = \left(\int_0^{\infty}U_-(\d x)\overline \Pi(x+ y)\right) \d y\,.$$
According to Vigon's {\it \'equation amicale invers\'ee}
(see \cite{Vi}), the right-hand side can be expressed as $\overline \mu_+(y)\d y$, 
where $\overline \mu_+$ denotes the tail of the L\'evy measure of the ladder height subordinator $H_+$. Hence we also have
$$\rho_2(\d y)= \frac{1}{E(H_+(1))} \left( a_+\delta_0(\d y) + \overline \mu_+(y) \d y\right) \,,$$
which is the classical limit distribution for the overshoot (i.e. the residual lifetime in the renewal process constructed from the subordinator $H_+$).

  More precisely, it belongs to the folklore of fluctuation theory  that when the ladder height of  a random walk or a L\'evy process has a finite mean, then the pair formed by the undershoot and the overshoot across a large level $z$ converges weakly as $z\to\infty$;  see in particular Theorem 3
 in \cite{KPR}.
 For every path $\omega\in \Omega$,  let 
 $$\hat T(\omega)= T(\hat \omega)=\inf\{t\geq 0: \omega(t)\in (0,\infty)\}$$ denote the first entrance time to the positive half-line. We now provide a formal statement of the convergence  alluded above which stresses the role of the measure $\rho$ defined in \eqref{eq5}.
 
 \begin{lemma}\label{L3} The probability measures on $[0,\infty)^2$
 $$ \hat P_z(\xi_{T-}\in \d x, -\xi_T\in \d y) 
 = P_{-z}(-\xi_{\hat T-}\in \d x, \xi_{\hat T}\in \d y)$$
 converge to $\rho$ as $z\to \infty$
 in the sense of weak convergence of probability measures.
 \end{lemma}
 \proof A by-product of the quintuple identity for first passage times  for L\'evy processes  (see Doney and Kyprianou \cite{DK} and references therein for further results in this area) is that for $x,y>0$
 $$P_{-z}(-\xi_{\hat T-}\in \d x, \xi_{\hat T}\in \d y)
 = \int_0^z U_+(z-\d v) {\bf 1}_{\{x>v\}}U_-(\d x-v)\Pi(\d x+y )\,.$$
 Roughly speaking, the renewal theorem implies that $U_+(z-\d v)$ converges as $z\to\infty$
 towards $\d v/E(H_+(1))$, and this yields that  
 $$\lim_{z\to \infty}P_{-z}(-\xi_{\hat T-}\in \d x, \xi_{\hat T}\in \d y)
 = \frac{1}{E(H_+(1))} U_-(x)\Pi(x+\d y) \d x\,,\qquad\hbox{vaguely on }  (0,\infty)^2.
$$
This establishes the claim when the ladder height process has no drift, because the right-hand side then defines a probability measure on $(0,\infty)^2$. 
In the case when  $a_+>0$,  the same conclusion follows invoking further Lemma \ref{L1}(i) and 
 the Portemanteau theorem, as $\rho$ is a probability measure. \QED

 \subsection{Time-reversal identities}
 On our way to providing the probabilistic interpretation of the duality relation of Theorem \ref{T1}, we 
 need to introduce some notation for c\`adl\`ag paths indexed by the entire real line.
 
 We set $\overline \Omega=\D(\R,\R)$;  $\overline \omega$ will denote a generic path in $\overline \Omega$. We also set $\xi_t(\overline \omega) = \overline \omega(t)$ for every $t\in \R$, so $(\xi_t)_{t\in\R}$ is  the usual coordinate process. 
  It may sometimes be convenient to identify $\overline \Omega$ as the product space
 $\Omega\times \Omega$ via the canonical bijection $\overline \omega \to (\omega, \omega')$, 
 where 
 $$\omega(t) = -\omega(-t-)\ \hbox{and}\ \omega'(t) = \omega(t)\qquad \hbox{for every }t\geq 0\,.$$
 Equivalently, we have for $t\in\R$
 $$\overline \omega(t) = \left\{\begin{matrix}
 \omega'(t)\ \hbox{if }t\geq 0\,,\\
 -\omega(-t-)\  \hbox{if } t<0\,.
 \end{matrix} \right. $$
 
 We then introduce a probability measure ${\mathcal P}$ on $\overline \Omega$ by
 $${\mathcal P}(\d \overline \omega)= {\mathcal P}(\d \omega, \d \omega') = \int_{[0,\infty)^2}\rho(\d x, \d y) P\up_x(\d \omega)P_y(\d \omega')\,.$$
 Thus under ${\mathcal P}$, the pair $(-\xi_{0-},\xi_0)$ has the stationary distribution $\rho$ for the under and the over shoots, and conditionally on $(-\xi_{0-},\xi_0)=(x,y)$,
 the processes $(-\xi_{t-})_{t\geq 0}$ and 
 $(\xi_t)_{t\geq 0}$ are independent with laws $P\up_x$ and $P_y$, respectively.
 Note  that $\xi_t<0$ for $t<0$ and  $\lim_{t\to -\infty}\xi_t=-\infty$ ${\mathcal P}$-a.s.
 while  for large times $\xi$ oscillates or 
 $\lim_{t\to \infty}\xi_t=\infty$ according as the genuine L\'evy process oscillates or tends to $\infty$.
In this section, we will essentially work with the canonical process on nonnegative times, $(\xi_t)_{t\geq 0}$, which has thus the law $P_{\rho_2}=\int \rho_2(\d x)P_x$ under ${\mathcal P}$.

We then fix a level $z>0$ and let $\tau(z)=\inf\{t\in\R: \xi_t>z\}$
denote the first passage time above $z$. Note that $\tau(z)\in[0,\infty)$ ${\mathcal P}$-a.s.
and that $\tau(z)=0$ when $\xi_0>z$. In the latter case, 
the notation $-\xi_{\tau(z)-}$ means $-\xi_{0-}$ (which is a nonnegative random variable),
and $(\xi_t)_{0\leq t <\tau(z)}$ is the empty path.

\begin{theorem}\label{T2} The following assertions hold for each $z>0$ :

\noindent {\rm (i)} We have the identity
$${\mathcal P}(z-\xi_{\tau(z)-}\in \d x, \xi_{\tau(z)}-z\in \d y) = \rho(\d x, \d y)\,.$$

\noindent{\rm (ii)} Under ${\mathcal P}$, the process $(\xi_t)_{0\leq t < \tau(z)}$ and the variable $\xi_{\tau(z)}$ are conditionally independent given $\xi_{\tau(z)-}$.

\noindent{\rm (iii)} Under the conditional law ${\mathcal P}(\cdot \mid \xi_{\tau(z)-}=z-x)$, the process $(z-\xi_{(\tau(z)-t)-})_{0\leq t < \tau(z)}$
has the same law  as the process $(\xi_t)_{0\leq t <\ell(z)}$
under $P\up_x$, where
$\ell(z)
=\sup\{t\geq 0: \xi_t<z\}$ denotes the last exit-time  from $(-\infty,z)$. 

\end{theorem}

\proof  (i) This just reflects the fact that $\rho$ is the stationary distribution of the  under and over shoots, as it can be seen from Lemma \ref{L3}.

(ii) Consider a process $K=(K_t)_{t \geq 0}$
with nonnegative left-continuous paths, which is adapted to the natural filtration generated by $\xi$ (hence $K$ is predictable). Let also $f: \R\to \R_+$ be a continuous function with support 
in $(z,\infty)$. We write $\Delta_t= \xi_t-\xi_{t-}$ for the jump of $\xi$ at time $t$ and ${\mathcal E}$ for the mathematical expectation under ${\mathcal P}$. By an application of the compensation formula to the Poisson point process of the jumps of a L\'evy process at the second line below, we have
\begin{eqnarray*}
{\mathcal E}\left(K_{\tau(z)}f(\xi_{\tau(z)})\right) &=&
{\mathcal E}\left(\sum_{t\geq 0}  {\bf 1}_{\{\sup_{0\leq s < t}\xi_s\leq z\}} 
K_t f(\Delta_t+\xi_{t-})\right)\\
&=&{\mathcal E}\left(\int_0^{\infty}\d t  {\bf 1}_{\{\sup_{0\leq s < t}\xi_s\leq z\}} K_t\left(\int \Pi(\d x)f(x +\xi_{t-}) \right)\right)\,.
\end{eqnarray*}
So if we define 
$$\varphi(y)=\frac{1}{\overline \Pi(z-y)}\int \Pi(\d x)f(x +y)  \qquad \hbox{for every }y\leq z$$
(we stress that $\varphi(z)=0$ when $\Pi((0,\infty))=\infty$, 
since our assumptions on $f$ ensure that $\int \Pi(\d x)f(x +z)<\infty$), 
we obtain
\begin{eqnarray*}
{\mathcal E}\left(K_{\tau(z)}f(\xi_{\tau(z)})\right) 
&=&{\mathcal E}\left(\int_0^{\infty}\d t  {\bf 1}_{\{\sup_{0\leq s < t}\xi_s\leq z\}} K_t \varphi(\xi_{t-})
\overline \Pi(z-\xi_{t-})\right)
\\
&=&
{\mathcal E}\left(\sum_{t\geq 0}  {\bf 1}_{\{\sup_{0\leq s < t}\xi_s\leq z\}} 
K_t \varphi(\xi_{t-}) {\bf 1}_{\{\Delta_t+\xi_{t-}>z\}}\right)\\
&=& {\mathcal E}\left(K_{\tau(z)}\varphi(\xi_{\tau(z)-})\right)\,,
\end{eqnarray*}
where we again applied the compensation formula for Poisson point processes at the second line above. This establishes the claim of conditional independence.

(iii)  To prove that under $P\up_{\rho_1}=\int\rho_1(\d x) P\up_x$ the process time-reversed at its last passage time below level $z$, $(\xi_{(\ell(z)-t)-}: 0\leq t<\ell(z))$, is sub-Markovian with semigroup $\hat p_t$ we shall apply Nagasawa's Theorem 3.5 in \cite{Na}. There are six conditions to be checked to ensure the validity of Nagasawa's result. To be rigorous and comply with the notation in \cite{Na} we put $E=(0,\infty)$ when $a_{+}=0$ and $E=[0,\infty)$ when $a_{+}>0$ ( note that $\rho_1$ has mass at zero and $\xi$ under $P\up_{\rho_{1}}$ can start from zero; however as mentioned before Lemma \ref{L1}, $P\up_{0}(\xi_{t}>0)=1$ for each $t>0$). We proceed by stating and demonstrating these conditions :

(1) Condition A 3.1 on p. 188 in \cite{Na} requires that the semigroups $p\up$ and $\hat{p}$ are in duality w.r.t. $\int_{0}^{\infty}U\up(x,dy)\rho_{1}(\d x)$. This holds due to Theorem \ref{T1} and $\rho_{1}(\d x)=m(\d x)/E(H_{+}(1))$.

(2) Condition A 3.2 and condition (ii) in the footnote on p. 190 in \cite{Na} are satisfied since our process is started from the probability measure $\rho_{1}$ and this ensures the finiteness of any quantities of the type $P\up_{\rho_{1}}(\xi_{(\ell(z)-t)-}\in.)$ and $\int_{0}^{\infty}e^{-\alpha t}P\up_{\rho_{1}}(\xi_{(\ell(z)-t)-}\in\cdot )\d t$, for $\alpha>0$.

(3) Condition A 3.3 asks for right-continuity of $\hat{P}^{T}_{t}f(x)=\hat{E}_{x}(f(\xi_{t}), T>t)$ in $t$ and a.s. right-continuity in $t$ of $\int_{0}^{\infty}e^{-\alpha s}\hat{P}^{T}_{s}f(\xi_{(\ell(z)-t)-})\d s$ under any of the measures $P\up_{a}$ for $a\in E$,  $x\in E$ and $\alpha>0$, where $f$ is a continuous function with a compact support in $E$. These are satisfied since the canonical process is defined in $\Omega$.

(4) Condition (i) in the footnote on p. 190 in \cite{Na} can be read off as the  existence 
of the left-limit $\xi_{\ell(z)-}$ on the event $\{0<\ell(z)<\infty\}$, $P_{\rho_{1}}$-a.s. This plainly holds.

(5) Condition (iii) is similar to A 3.3. We need to check that $\hat{P}^{T}_{t}f(x)$ is continuous for $x\in E$ for any continuous function $f$ with a compact support in $E$. However the statement is clear from the continuity of $\hat{P}_{x}(T>t)$ in $x$, for $x>0$ and any fixed $t$, independently of $a_{+}=0$ or not, and the fact that when $a_{+}>0$, $\hat{T}$ is continuous at $0$ since $0$ is regular for $(0,\infty).$

  Thus, we have checked that $(\xi_{(\ell(z)-t)-}: 0\leq t<\ell(z))$ is sub-Markovian with semigroup $\hat p_t$. Note that the latter thus fulfills the strong Markov property and 
denote its initial distribution by
$$\nu_z(\d x)=P\up_{\rho_1}(\xi_{\ell(z)-}\in \d x, \ell(z)>0)\,, \qquad\hbox{for }x\in[0,z].$$
Pick an arbitrary $z'>z$ and apply the strong Markov property to the time-reversed path $(\xi_{(\ell(z')-t)-}: 0\leq t<\ell(z'))$ at its first-passage time below level $z$, that is
$\ell(z')-\ell(z)$. This yields
$$\nu_z(\d x)=\int_0^{\infty}\nu_{z'}(\d y) \hat P_{y-z}(z+\xi_T\in \d x)\,, \qquad\hbox{for }x\in[0,z].$$
Letting $z'\to\infty$, we deduce from Lemma \ref{L3} that $\nu_z(\d x) = \rho_2(z-\d x)$
on $[0,z]$.

We conclude that under $P\up_{\rho_1}$,  the  time-reversed  process $(\xi_{(\ell(z)-t)-})_{ 0\leq t<\ell(z)}$ has the same law as $(\xi_t)_{0\leq t < T}$ under $\hat P_{\nu_z}$.
Equivalently, under $P_{\rho_2}$, the process $(z-\xi_{(\tau(z)-t)-})_{0\leq t < \tau(z)}$
has the same law as the process $(\xi_t)_{0\leq t <\ell(z)}$
under $P\up_{\rho_1}$, which is our statement. \QED

In the special case when the drift coefficient $a_+$ of the ladder height subordinator $H_+$
is strictly positive, recall from Lemma \ref{L1} that the probability under $\hat P_x$ that the first-exit from $(0,\infty)$ occurs continuously equals $a_+u_+(x)>0$. 
Similarly, it is easy to see that for any $x>0$, the probability under $P\up_0$ that the last-exit from $(0,x)$ occurs continuously is strictly positive. Hence we immediately deduce
from Theorem \ref{T2} the following identity which can also be seen from Theorem 4 of Chaumont \cite{Ch}.

\begin{corollary}\label{C1} Suppose $a_+>0$ and fix $x>0$. 
The law of $(\xi_t: 0\leq t<\ell(x))$ under $P\up_{0}(\cdot \mid \hbox
{ $\xi$ is continuous at $\ell(x)$})$ is
that of $(\xi_{(T-t)-}: 0\leq t<T)$ under $\hat P_{x}(\cdot \mid \hbox{ $\xi$ is continuous at }T)$.
\end{corollary}
{\bf Remark.} Note that these conditionings are trivial when $\xi$ has no positive jumps under $P$, and hence Corollary \ref{C1} encompasses the extension of Williams' time-reversal stated as Theorem VII.18 in \cite{Be}.

In the same vein, we recover (and slightly extend) a result due to Duquesne;
see Theorems 4.1 and 4.2 in \cite{Du}.

\begin{corollary}\label{C2}
Under $P\up_{0}$ the law of $(\xi_{\ell(x)-}-\xi_{(\ell(x)-t)-})_{t<\ell(x)}$ is the law of $(\xi_{t})_{t<g(\hat{T}(x))}$ under $P$, where $\hat{T}(x)=\inf\{t\geq0:\xi_{t}>x\}$ and $g(\hat{T}(x))=\sup\{s<\hat{T}(x):\xi_{s}\vee \xi_{s-}=\sup_{t<\hat{T}(x)}\xi_{t}\}$.
\end{corollary}
{\bf Remark.} Note that Corollary \ref{C2} is intuitively obvious when we deal with random walks because the law $P\up_0$ equals the law under $P$ of the reversed excursions away from the maximum.

\proof Recall that we assume that $E(H_{+}(1))<\infty$ and denote by
$I=\inf_{t\geq0} \xi_{t}\text{ and }g_{I}=\sup_{s\geq0}\{\xi_{s}\wedge\xi_{s-}=I\}$. We deduce from Theorem 25 in Chapter 8 in \cite{Do} that under $P\up_{\rho_{1}}$ the process $(\xi_{g_{I}+t}-\xi_{g_{I}})_{t\geq 0}$ has the law of $(\xi_{t})_{t\geq0}$ under $P\up_0$ and it is independent of $(\xi_{t})_{t\leq g_{I}}$. This combined with Theorem \ref{T2} yields that for every bounded measurable functional $F$ on $\Omega$
\begin{eqnarray}\label{M1}
E\up_{\rho_1}\left[F\big((\xi_{\ell(x)-}-\xi_{(\ell(x)-t)-})_{t<\ell(x)-g_I}\big)\right]
&=&E\up_0\left[F\big((\xi_{\ell(x-I)-}-\xi_{(\ell(x-I)-t)-})_{t<\ell(x-I)}\big)\right] \nonumber \\
&=&E\left[F\big((\xi_{t})_{t<g(T(x-\rho_{2}))}\big)\right],\end{eqnarray}
where $g(T(x-\rho_{2}))$ is the time when the last maximum is attained before $\xi$ goes beyond $x-\rho_{2}$, where $\rho_{2}$ is the stationary overshoot and is independent of $\xi$ under $P$ and $I$ in the second term has the distribution of $P\up_{\rho_1}(I\in dy)$ and is independent of $\xi$ under $P\up_0$. Next note that Theorem \ref{T2} implies as well that 
\[P\up_{\rho_1}(I\in dy)=\lim_{z\to\infty}\hat{P}_{\rho_1}(\inf_{s<\ell(z)}\xi_{s}\in dy)=\lim_{z\to\infty}P_{\rho_{2}}(\inf_{s<\tau(z)}z-\xi_{\tau(z-s)-}\in dy)=\rho_{2}(dy).\]
 The last identity can be deduced as in Lemma \ref{T3} or recovered from Theorem 3 in \cite{KPR}. Then \eqref{M1} translates easily to $A*\rho_{2}(x)=B*\rho_{2}(x)$ with $B(y)=EF\big((\xi_{t})_{t<g(T(y))}\big)$ and $A(y)=E\up_0F\big((\xi_{\ell(y)-}-\xi_{(\ell(y)-t)-})_{t<\ell(y)}\big)$ and we get that $A=B$ on $(0,\infty)$. Since this holds for any measurable functional $F$ we conclude the proof.
\QED
 \end{section}
  
 \begin{section} {Applications to weak limit theorems}
 \subsection{Starting a L\'evy process from $-\infty$}

 We first observe that the probability measure ${\mathcal P}$ that has been introduced in Section 2.4 fulfills a remarkable spatial stationarity property,
 which follows easily from Theorem \ref{T2} and the strong Markov property.
  
\begin{corollary}\label{C3} For any $x\in\R$, let $\tau(x)=\inf\{t\in\R: \xi_t>x\}$
denote the first passage time of $\xi$ above the level $x$. Under ${\mathcal P}$, the processes $(\xi_{\tau(x)+t})_{t\in\R}$ and $(x+\xi_t)_{t\in\R}$ have the same distribution.
\end{corollary}
 
\proof  Let   ${\mathcal P}_z$ be the law of  $(z+\xi_t)_{t\in\R}$
for an arbitrary $z\in\R$. An application of the strong Markov property combined with
Theorem \ref{T2} shows the following. Let us work under  ${\mathcal P}_z$ for an arbitrary $z<0$, and recall that $\hat T$ denotes the first entrance time into $(0,\infty)$. Then
the pair $(-\xi_{\hat T-}, \xi_{\hat T})$ has the law $\rho$, and conditionally on
$(-\xi_{\hat T-}, \xi_{\hat T})=(x,y)$, the processes 
$(-\xi_{(\hat T-t)-})_{0\leq t < \hat T}$ and $(\xi_{t+\hat T})_{t\geq 0}$ are independent.
Further the former has the same law as $(\xi_t)_{0\leq t < \ell(-z)}$ under $P\up_x$
while the latter has the law $P_y$. As the last-exit time from $(-\infty,-z)$, $\ell(-z)$, tends to infinity as $z\to -\infty$, we see that
the law of the shifted process $(\xi_{t+\hat T})_{t\in\R}$ under
${\mathcal P}_z$ converges weakly to that of $(\xi_t)_{t\in\R}$ under ${\mathcal P}$. 
We can now complete the proof by an easy argument based on replacing $z$ by $z+x$. \QED

 The proof of Corollary \ref{C3} suggests that  the limit theorem for the under and over shoots should have an extension to paths. In this direction, for every $\omega\in \D([0,\infty),\R)$,
we denote  by  $\vartheta (\omega)$ the path indexed by the entire real line which is obtained by shifting $\omega$ at the time $\hat T(\omega)=\hat T$ of its first entrance in $(0,\infty)$, that is  $\vartheta (\omega)=(\omega'(s))_{s\in\R}$  with
  $$ \omega'(s)=\left\{ \begin{matrix} \omega(\hat T+s) \hbox{ for every }s\geq -\hat T\,,\\
 -\infty \hbox{ otherwise}\,.
 \end{matrix}\right.
 $$
 
We now state the main result of this section. 
 \begin{theorem} \label{T3}  Fix any $b\in\R$. The law of $(\xi_t\circ \vartheta )_{t\geq b}$ under $P_x$
 converges weakly on $\D([b,\infty),\R)$ as $x\to -\infty$ towards the law
 of $(\xi_t)_{t\geq b}$ under ${\mathcal P}$.
 \end{theorem}
 
 We shall derive Theorem \ref{T3} from Corollary \ref{C3}  using a coupling argument which requires distinguishing whether the drift coefficient
 $a_+$ of the ladder height subordinator $H_+$ is zero or strictly positive.
 The first case is easier and relies on the following standard construction.
 
 \begin{lemma}\label{L4} Fix $\varepsilon >0$. We can construct  on some probability space
 a random variable $\gamma$ with values in $[0,\varepsilon]$, an a.s. finite random time $\tau$
 and a pair of processes $(\xi'_t)_{t\geq 0}$ and $(\xi''_t)_{t\geq 0}$ that fulfill the following requirements :
 
 \noindent {\rm (i)} $\xi'$ has the law $P$ and $\xi''$ has the law $P_{\rho_2}$,
 
  \noindent {\rm (ii)} $\xi''_t=\xi'_t+\gamma$ for all $t\geq \tau$.
 
 \end{lemma}
 
  \proof We start from a pair $(\tilde \xi', \xi'')$ with law $P\otimes P_{\rho_2}$, i.e. 
 $\tilde \xi' $ and $\xi''$ are independent with respective laws $P$ and $P_{\rho_2}$.
 Then $\xi''-\tilde \xi'$ is a symmetric L\'evy process with initial law $\rho_2$. Recall from Lemma \ref{L2}  that it is centered and hence recurrent by the test of Chung and Fuchs (cf. Exercise I.10 in \cite{Be}). We set $\tau=\inf\{t\geq 0: \xi''_t-\tilde \xi'_t \in[0,\varepsilon]\}$ and $\gamma= \xi''_{\tau}-\tilde \xi'_{\tau}$, so $\tau$ is an a.s finite
 stopping time and $\gamma$ a random variable in $[0,\varepsilon]$. By the strong Markov property, the process 
 $$\xi'_t=\left\{\begin{matrix}\tilde \xi'_t \ \hbox{ if }t \leq \tau\\
\xi''_t- \gamma\  \hbox{ if }t> \tau
 \end{matrix}\right.
 $$
 has  the law $P$. \QED
 
When  the drift coefficient  of the ladder height subordinator $H_+$ is strictly positive,  $a_+>0$,
we need a stronger coupling. 
 
 \begin{lemma}\label{L5} Assume $a_+>0$.  We can construct  on some probability space
 a pair of a.s. finite random times $\tau'$ and $\tau''$
 and a pair of processes $(\xi'_t)_{t\geq 0}$ and $(\xi''_t)_{t\geq 0}$ that fulfill the following requirements :
 
 \noindent {\rm (i)} $\xi'$ has the law $P$ and $\xi''$ has the law $P_{\rho_2}$,
 
  \noindent {\rm (ii)} $\xi'_{\tau'+t}=\xi''_{\tau''+t}$ for all $t\geq 0$.
 
 \end{lemma}
 
 \proof  We start again from a pair $(\tilde \xi', \xi'')$ with law $P\otimes P_{\rho_2}$.
 We define  the passage times
 $$\sigma'_1=\inf\{t\geq 0: \tilde \xi'_t\geq \xi''_0\}\ ,\ \sigma''_1=\inf\{t\geq 0: \xi''_t\geq \tilde \xi'(\sigma'_1)\}\,$$
 and then recursively
  $$\sigma'_{k+1}=\inf\{t\geq 0: \tilde \xi'_t\geq \xi''(\sigma''_k)\}\ ,\ \sigma''_{k+1}=\inf\{t\geq 0: \xi''_t\geq \tilde \xi'(\sigma'_{k+1})\}\,.$$
  
We claim that a.s., these non-decreasing sequences remain constant after a finite number of steps. 
Taking this assertion for granted, the construction of the coupling is immediate as it suffices to set
$\tau'= \sigma'_{\infty}$, $\tau''=\sigma''_{\infty}$, and 
$$\xi'_t=\left\{\begin{matrix}\tilde \xi'_t \ \hbox{ if }t \leq \tau'\\
\xi''_{\tau''+t-\tau'}\  \hbox{ if }t> \tau'\,.
 \end{matrix}\right.
 $$
To complete the proof, it suffices to recall from Lemma \ref{L1} (i) that, since 
the ladder height subordinator $H_+$ has a strictly positive drift coefficient and a finite mean, the probability that $H_+$ hits some fixed 
point $x\in[0,\infty)$ is bounded from below by a strictly positive constant. 
Hence the number of steps alluded above is stochastically bounded  by a geometric variable.
\QED

 We may now tackle the proof of Theorem \ref{T3}.
 
\proof  Recall the notation ${\mathcal P}_x$ for the law of $x+\xi$ under ${\mathcal P}$

1. Suppose first $a_+>0$ and fix $\varepsilon>0$ arbitrarily small. According to Lemma \ref{L5}, we can construct two process $(\xi'_t)_{t\geq 0}$ and $(\xi''_s)_{s\in\R}$ with
respective laws $P$ and ${\mathcal P}$
and two a.s. finite random times $\tau'$ and $\tau''$ such that 
$\xi'_{\tau'+t}=\xi''_{\tau''+t}$ for all $t\geq 0$. Provided that $x$ is chosen sufficiently large, the probability of the event that
$$\sup_{0\leq t \leq \tau'}\xi'_t \leq x/2 \ , \  \sup_{ t \leq \tau''}\xi''_t \leq x/2 
\hbox{ and } 
\inf\{t\geq 0: \xi'_{t+\tau'}-\xi'_{\tau'}>x/2\} > -b$$
is at least $1-\varepsilon$. 

Therefore if we set $\tilde \xi'_t=\xi'_t-x$ and $\tilde \xi''_t=\xi''_t-x$, then the processes
$(\tilde \xi'_t)_{t\geq 0}$ and $(\tilde \xi''_s)_{s\in\R}$ have the
 law $P_{-x}$ and ${\mathcal P}_{-x}$, respectively. Further the probability
 that the paths obtained by shifting $\tilde \xi'$ and $\tilde \xi''$ at their first entrance time into $(0,\infty)$ coincide on $[b,\infty)$ is bounded from below by $1-\varepsilon$. 
This entails the statement as we know from Corollary \ref{C3} that the path obtained by shifting $\tilde \xi''$ has the law ${\mathcal P}$.

2. Suppose now that $a_+=0$ and fix $\eta>0$ arbitrarily small. Then the stationary distribution $\rho_2$ of the overshoot has no atom at $0$, and we may pick $\varepsilon>0$ sufficiently small so that $\rho_2([0,\varepsilon])<\eta/2$. 
According to Lemma \ref{L4}, we can construct two process $(\xi'_t)_{t\geq 0}$ and $(\xi''_s)_{s\in\R}$ with
respective laws $P$ and ${\mathcal P}$, 
an a.s. finite random time $\tau$ and a random variable $\gamma\in [0,\varepsilon]$ such that $\xi'_{t}+\gamma=\xi''_{t}$ for all $t\geq \tau$. 

For every $x\geq 0$, consider the first entrance times
$$\tau'(x)=\inf\{t\geq 0: \xi'_t>x\}\ \hbox{and}\ \tau''(x)=\inf\{t\geq 0: \xi''_t>x\}\,.$$
Since $\xi'_t\leq \xi''_t \leq \xi'_t+\varepsilon$ for all $t>\tau$,  
the probability of the event 
$$\{\tau<\tau'(x)\wedge \tau''(x), \tau'(x)\neq \tau''(x)\}$$
is bounded from above by the probability that the overshoot $\xi''_{\tau''(x)}-x$ does not exceed $\varepsilon$, and thus by $\eta/2$. Further, provided that $x$ is chosen sufficiently large, the probability of the event $\{\tau-b<\tau'(x)\wedge \tau''(x)\}$
is at least $1-\eta/2$. 

Now set $\tilde \xi'_t=\xi'_t-x$ and $\tilde \xi''_t=\xi''_t-x$; the processes
$(\tilde \xi'_t)_{t\geq 0}$ and $(\tilde \xi''_s)_{s\in\R}$ have thus the
 law $P_{-x}$ and ${\mathcal P}_{-x}$, respectively. It follows from above that the probability that the paths obtained by shifting $\tilde \xi'$ and $\tilde \xi''$ at their first entrance time into $(0,\infty)$ remain parallel  on the time interval $[b,\infty)$ with a distance at most $\varepsilon$,  is bounded from below by $1-\eta$. 
This entails the statement as we know from Corollary \ref{C3} that the path obtained by shifting $\tilde \xi''$ has the law ${\mathcal P}$. \QED

 \subsection{A Lamperti-type representation for self-similar Markov processes
 entering from $0+$}
 We now conclude this work with an application to the class of Markov processes in $(0,\infty)$ that enjoy the scaling property. That is, we consider a Markov process
 $X=(X_t)_{t\geq 0}$ with values in $(0,\infty)$ and write ${\P_x}$ for its law started from 
 $X_0=x>0$. We shall always assume that the process is conservative, i.e. there is no cemetery state. We suppose that the self-similarity property
 $$\hbox{the distribution of $(cX_{t/c})_{t\geq 0}$ under $\P_x$ is $\P_{cx}$}$$
 holds for every $c,x>0$. Lamperti \cite{La} has studied in depth this class of processes
 which are nowadays called positive self-similar Markov processes (in short pssMp),
 and pointed at a fundamental connection with real valued L\'evy processes that can be described as follows in the framework of this paper.
 
 We work under the probability measure $P_y$ for which $(\xi_t)_{t\geq 0}$
 is a L\'evy process started from $y\in\R$. We drop for a moment the assumption
 that the ladder height has a finite expectation, and just suppose that the L\'evy process does not tend to $-\infty$. We introduce a time-change $\gamma(t)$  for every $t\geq 0$ as the inverse of the exponential functional, that is
 $$\int_0^{\gamma(t)} \e^{\xi_s} \d s\,=t\,.$$
 Then the process $X_t=\exp( \xi_{\gamma(t)})$ is a pssMp started from $x=\e^y$,
 and any (conservative) pssMp can be constructed by this procedure. 
 
 The question of whether a pssMp can enter from the boundary point $0+$,
 that is if $\P_x$ admits a non-degenerate weak limit $\P_{0+}$ as $x\to 0+$
 was raised by Lamperti. Bertoin and Yor \cite{BY} provided a positive 
 answer when the underlying L\'evy process possesses a positive and finite first moment.
 Recall that this implies that the mean ladder height $E(H_+(1))$ is finite; it is further easy to show that a pssMp cannot enter from $0+$ when $E(H_+(1))=\infty$.
 Caballero and Chaumont \cite{CC} obtained an explicit necessary and sufficient condition; basically they proved that a pssMp can enter from $0+$ if and only if $E(H_+(1))<\infty$ and some very mild technical condition holds. That this technical condition
 is automatically fulfilled when the mean ladder height is finite has been proved recently by Chaumont {\it et al.} \cite{CKPR}, so the definitive simple characterization is that
 a pssMp can enter from $0+$ if and only if  $E(H_+(1))<\infty$. 
>From now on this assumption is thus again enforced. 
 
 We point at a simple and direct construction of the law $\P_{0+}$ in terms of the spatially homogeneous law ${\mathcal P}$ and the canonical process $(\xi_t)_{t\in\R}$ indexed by the whole real line, which is somehow hidden in the approach by Caballero and Chaumont \cite{CC}.  The intuition stems from  the observation that Lamperti's transformation can be re-expressed in terms of the shifted path $\xi\circ\vartheta$: if we write $\sigma: [0,\infty)\to \R$ for the inverse of 
  the functional
  $$t\to \int_{-\infty}^t\exp(\xi_s\circ\vartheta) \d s\,, \qquad t\in \R\,,$$
  then for every $y\in\R$, under $P_y$ the time-changed process
  $(\exp(\xi_{\sigma(t)}\circ \vartheta))_{t\geq 0}$ has the law $\P_x$ with $x=\e^y$. 
  Since we know from Theorem \ref{T3} that the law of the shifted path $\xi\circ\vartheta$
  under $P_y$ converges weakly to ${\mathcal P}$, we arrive naturally at the following.

 \begin{corollary}\label{C4} {\rm (i)} The exponential functional 
 $$I(t)=\int_{-\infty}^t \e^{\xi_s}\d s$$
 is finite for all $t\in\R$ and $I(\infty)=\infty$, ${\mathcal P}$-a.s. 
 
 \noindent{\rm (ii)} Introduce the time-change $\sigma(t)$ for $t> 0$ by
  $$\int_{-\infty}^{\sigma(t)} \e^{\xi_s} \d s\,=t\,.$$
  Then under ${\mathcal P}$,  the process $(X_t)_{t>0}$ given by
  $$X_t=\exp( \xi_{\sigma(t)})$$
  is a pssMp started from the entrance boundary $0+$,
  in the sense that its law is the weak limit of $\P_x$ as $x\to 0+$.
  \end{corollary}

 \proof (i) By Theorem \ref{T1}(ii), we have
 $${\mathcal E}(I(0))=\int_0^{\infty} \rho_1(\d x) \int_0^{\infty}U\up(x,\d y)\e^{-y}
 = \frac{1}{E(H_+(1))} \int_0^{\infty}\d y \e^{-y} U_-(y)\,.$$
 Because a renewal function is always sub-additive, the right-hand side is finite,
 which implies our claim. 
 
 (ii) We now know that $X_t=\exp(\xi_{\sigma(t)})$ is well-defined. Recall that
 $\tau(y)=\inf\{s\in\R: \xi_s>y\}$ denotes the first entrance time in $(y,\infty)$, so
 $I(\tau(y))$ corresponds to the first passage time of $X$ above the level $x=\e^y$. It follows readily from Corollary  \ref{C3} and Lamperti's transformation that  conditionally on $X_{I(\tau(y))}=z$,
 the shifted process $(X_{t+I(\tau(y))})_{t\geq 0}$ has the law $\P_{z}$. Further the distribution of $X_{I(\tau(y))}$, say
$\mu_x$, is the image of  $\rho_2$ by the map $w\to x\e^{w}$. Plainly $\mu_x$ converges to $\delta_0$ as $x\to 0+$ (i.e. $y\to -\infty$), while $\tau(y)=\sigma(I(\tau(y))\to -\infty$.
We conclude that $(X_t)_{t>0}$ has the law $\P_{0+}$. \QED

We stress that the argument for establishing Corollary \ref{C4} has little to do with self-similarity or exponential functions. It would apply just as well to construct the Markov process entering from the boundary $0+$ in the following more general situation 
(which mirrors Feller's construction of one-dimensional diffusions as time-space transforms of Brownian motions). More precisely,  consider a measurable locally bounded  function
$f: \R\to (0,\infty)$ with $\int_{-\infty}^0 \d y |y|f(y)<\infty$ 
and $\int_0^{\infty}\d y f(y) = \infty$, and $g:\R\to (0,\infty)$ a continuous strictly increasing function with $\lim_{z\to-\infty}g(z)=0$. 
Then the functional
$$I(t)=\int_{-\infty} ^t f(\xi_s)\d s$$
is finite for all $t\in\R$ and $\lim_{t\to\infty}I(t)=\infty$, ${\mathcal P}$-a.s.
Writing $\sigma$ for the inverse functional of $I$, the process
$(X_t)_{t>0}$ defined by $X_t=g( \xi_{\sigma(t)})$ then enters from $0+$,
and is Markovian in $(0,\infty)$ with an infinitesimal generator defined by an obvious transformation of that of the L\'evy process. 
 \end{section}
 
 \vskip1cm \noindent {\bf
Acknowledgment :} This work has been undertaken while M. S. was on a post-doctoral 
position at the Laboratoire de Probabilit\'es et Mod\`eles Al\'eatoires, thanks to
a financial support from the University Pierre-et-Marie Curie.

  \end{document}